\documentclass[final]{siamltex}

\usepackage{color}
\usepackage{a4,amsmath,amssymb,amscd,latexsym} 
\usepackage{bbm} 
\usepackage{epsfig,color}
\usepackage{amsfonts}
\usepackage{mathrsfs}

\usepackage{overpic}

\usepackage{paralist}

\usepackage{graphicx}

\usepackage{trfsigns}

\usepackage{stmaryrd}

\def\hess#1{\mbox{\sl Hess}_{\! #1}\, }
\def\grad{\mbox{grad}}
\newcommand{\LL}{{\mathscr{L}}}

\def\bei#1{\,\rule[-4mm]{.1mm}{9mm}_{\hspace{.6mm} #1 }}

\newcommand{\ACS}{2pt}
\newcommand{\mat}[4]{{\arraycolsep\ACS
\left#1\begin{array}{@{}*{#2}{c}@{}}#4\end{array}\right#3}}

\def\cN{{\cal N}}

\def\exp{\mathrm{exp}}

\newcommand{\R}{{\mathbbm{R}}} 

\makeatletter
\def\moverlay{\mathpalette\mov@rlay}
\def\mov@rlay#1#2{\leavevmode\vtop{%
   \baselineskip\z@skip \lineskiplimit-\maxdimen
   \ialign{\hfil$\m@th#1##$\hfil\cr#2\crcr}}}
\newcommand{\charfusion}[3][\mathord]{
    #1{\ifx#1\mathop\vphantom{#2}\fi
        \mathpalette\mov@rlay{#2\cr#3}
      }
    \ifx#1\mathop\expandafter\displaylimits\fi}
\makeatother

\newcommand{\cupdot}{\charfusion[\mathbin]{\cup}{\cdot}}

\def\scp#1{\left\langle #1\right\rangle}

\def\CROP#1{}

\newtheorem{remark}{Remark}

\title{Towards a Lagrange--Newton Approach for PDE constrained Shape Optimization
}

\author{Volker H.~Schulz, Martin Siebenborn, and Kathrin Welker\thanks{University of Trier, Department of Mathematics, 54296 Trier, Germany ({\tt volker.schulz@uni-trier.de, siebenborn@uni-trier.de, welker@uni-trier.de}).}
        }

\begin{document}

\maketitle

\begin{abstract}
The novel Riemannian view on shape optimization developed in \cite{VHS-shape-Riemann} is extended to a Lagrange--Newton approach for PDE constrained shape optimization problems. The extension is based on optimization on Riemannian vector space bundles and exemplified for a simple numerical example.
\end{abstract}



\pagestyle{myheadings}
\thispagestyle{plain}

\section{Introduction}
Shape optimization problems arise frequently in technological processes, which are modeled in the form of partial differential equations as in \cite{Arian-1995,ARTASAN-1996,Arian-Vatsa-1998,ESSI-2009,GISS2012,ComFluid2011,SS-2009,CoCy2010}. In many practical circumstances, the shape under investigation is parameterized by finitely many parameters, which on the one hand allows the application of standard optimization approaches, but on the other hand limits the space of reachable shapes unnecessarily. Shape calculus, which has been the subject of several monographs \cite{Delfour-Zolesio-2001,Mohammadi-2001,Sokolowski-1992} presents a way out of that dilemma. However, so far it is mainly applied in the form of gradient descent methods, which can be shown to converge. The major difference between shape optimization and the standard PDE constrained optimization framework is the lack of the linear space structure in shape spaces. If one cannot use a linear space structure, then the next best structure is the Riemannian manifold structure as discussed for shape spaces in \cite{BHM-2011-unpublished,BHM-2011,Michor06anoverview,MM-2006,MM-2005}. The publication \cite{VHS-shape-Riemann} makes a link between shape calculus and shape manifolds and thus enables the usage of optimization techniques on manifolds in the context of shape optimization. 

PDE constrained shape optimization however, is confronted with function spaces defined on varying domains. The current paper presents a vector bundle framework based on the Riemannian framework established in \cite{VHS-shape-Riemann}, which enables the discussion of Lagrange--Newton methods within the shape calculus framework for PDE constrained shape optimization.

The paper first presents the novel Riemannian vector bundle framework on section \ref{sec2}, discusses this approach for a specific academic example in section \ref{sec3} and presents numerical results in section \ref{sec4} .

\section{Constrained Riemannian shape optimization}\label{sec2}
The typical set-up of an equality constrained optimization problem is
\begin{align*}
\min\limits_{y,u}\ & J(y,u)\, , \ J\colon Y\times U\to \R\\
\mbox{s.t.}\ & c(y,u)=0\, , \ c\colon Y\times U\to Z
\end{align*}
where $U,Y,Z$ are linear spaces and $c,J$ sufficiently smooth nonlinear functions \cite{optbook}. In some situations the constraint $c$ allows to apply the implicit function theorem in order the define a unique control to state mapping $y(u)$ and thus the constrained problem maybe reduced to an unconstrained one of the form
\[
\min\limits_u J(y(u),u).
\]
However, the constrained formulation is often computationally advantageous, because it allows the usage of pre--existing solver technology for the constraint and it is geared towards an efficient SAND (simultaneous analysis and design) or one--shot approach based on linear KKT systems. So far, shape optimization methods based on the shape calculus, have been mainly considered with the reduced black--box framework above via the implicit function theorem -- mainly because the set of all admissible shapes is typically not a linear space -- unlike the space $U$ above. The publication \cite{VHS-shape-Riemann} has developed a Riemannian framework for shape optimization in the reduced unconstrained paradigm, which enables Newton--like iteration techniques and convergence results. This publication aims at generalizing those results to the constrained perspective -- in particular  for the case that the constraint is of the form of a set of partial differential equations (PDE). 

Within that framework, the space $Y$ for the state variable is a linear (function) space depending explicitly on $u\in U$, e.g., $H^1(\Omega(u))$, where $\Omega(u)$ is the interior of a shape $u$. This line of thinking leads to vector bundles of function spaces as discussed in detail in \cite{Lang-2001}. Thus, we now consider a Riemannian manifold $(\cN,G)$ of class $C^q$ ($q\ge 0$), where $G$ is a smooth mapping assigning any point $p\in\cN$ an inner product $G_p(\cdot,\cdot)$ on the tangential bundle $T\cN$. 
For each $u\in\cN$, there is given a Hilbert space $H(u)$ such that the set
\[
E:=\left\{(H(u),u)\, |\, u\in\cN\right\}
\]
is the total space of a vector bundle $(E,\pi,\cN)$. In particular, there is a bundle--projection $\pi:E\to\cN$ and for an open covering $\left\{U_i\right\}$ of $\cN$ a local $C^q$ isomorphism
\[
\tau_i\colon\pi^{-1}(U_i)\to H_0 \times U_i 
\]
where $H_0$ is a Hilbert space. In particular, we have an isomorphism on each fiber
\[
\tau_i(u)\colon\pi^{-1}(x)=H(u)\to  H_0
\]
and for $u\in U_i\cap U_j$, the mapping $\tau_i(u)\circ\tau_j(u)^{-1}:H_0\to H_0$ is a linear isomorphism. The total space $E$ of the vector bundle $(E,\pi,\cN)$ is by itself a Riemannian manifold, where the tangential bundle $TE$ satisfies
\[
T_{(y,u)}E\cong  H(u)\times T_y\cN.
\]

In Riemannian geometry, tangential vectors are considered as first order differential operators acting on germs of scalar valued functions (e.g. \cite{conlon}). 
Such a differential operator will be notated by $h(J)(e)$, if $J\colon E\to\R$ is a differentiable function and $e\in E$.
We will have to deal with derivatives, where we will always use directional derivatives of scalar valued functions only, but notate them in the usual fashion. Let the derivative of $J$ at $e$ in direction $h\in TE$ be denoted by $DJ(e)h$. Then, we define in this setting
\[
DJ(e)h:=h(J)(e)\, , \ h\in TE .
\]
In particular, we denote
\begin{align*}
\frac{\partial}{\partial y}J(y,u)h_y&:=h_1(J)(y,u)\, ,\ h_1:=(h_y,0)\in TE\\ 
\frac{\partial}{\partial u}J(y,u)h_u&:=h_2(J)(y,u)\, ,\ h_2:=(0,h_u)\in TE 
\end{align*}
where $h_y\in H(u)$ and $h_u\in T_y\cN$, if $h_1, h_2\in T_{(y,u)}E$.

We consider now the following constrained optimization problem
\begin{align}\label{vb-opt1}
&\min \limits_{(y,u)\in E}\   J(y,u)\, , \ J\colon E\to \R\\
\label{vb-opt2}
&\hspace{0,4cm}\mbox{s.t.}\ a_u(y,p)=b_u(p)\, , \ \forall p\in H
\end{align}
where $a_u(.,.)$ is a bilinear form and $b_u(.)$ a linear form defined on the fiber $H$ which are $C^q$ with respect to $u$. The scalar valued function $J$ is assumed to be $C^q$. Intentionally, the weak formulation of the PDE is chosen for ease of presentation. Now, it will be necessary to define the Lagrangian $\LL$ in order to formulate the adjoint and design equation to the constrained optimization problem (\ref{vb-opt1}--\ref{vb-opt2}).

\begin{definition}\label{lagrangian}
We define the Lagrangian in the setting above for $(y,u,p)\in F$ as
\[
\LL(y,u,p):=J(y,u)+a_u(y,p)-b_u(p)
\]
where $F:=\left\{(H(u),u,H(u))\,|\,u\in\cN\right\}$ with $T_{(y,u,p)}F\cong H(u)\times T_y\cN\times H(u)$.
\end{definition}

Let $(\hat{y},\hat{u})\in E$ solves the optimization problem (\ref{vb-opt1}--\ref{vb-opt2}). Then, the (adjoint) variational problem which we get by differentiating $\LL$ with respect to $y$ is given by 
\begin{equation}\label{adjoint}
a_{\hat{u}}(z,p)=-\frac{\partial}{\partial y}J(\hat{y},\hat{u})z\, , \ \forall z\in H(\hat{u})
\end{equation} 
and the design problem which we get by differentiating $\LL$ with respect to $u$ is given by
\begin{equation}\label{design}
\frac{\partial}{\partial u}\bei{u=\hat{u}}\left[J(\hat{y},u)+a_{u}(\hat{y},\hat{p})-b_{u}(\hat{p})\right]w=0\, , \ \forall w\in T_{\hat{u}}\cN
\end{equation}
where $\hat{p}\in H$ solves (\ref{adjoint}). If we differentiate $\LL$ with respect to $p$, we get the state equation (\ref{vb-opt2}). These (KKT) conditions (\ref{vb-opt2}--\ref{design}) could be collected in the following condition:
\begin{equation}\label{KKT}
D\LL(\hat{y},\hat{u},\hat{p})h=0\, , \ \forall h\in T_{(y,u,p)}F.
\end{equation}

\begin{remark}
In a vector space setting, the existence of a solution $p\in H$ of the (adjoint) variational problem (\ref{adjoint}) is typically guaranteed by so--called constraint qualifications. From this point of view, here, the existence itself can be interpreted as formulation of a constraint qualification.
\end{remark}


By using a Riemannian metric $G$ on $T\cN$ and a smoothly varying scalar product $\scp{.,.}_u$ on the Hilbert space $H(u)$, we can envision $T_{(y,u,p)}F$ as a Hilbert space with a canonical scalar product 
\begin{equation}\label{scalarproduct}
\scp{
\mat(1){z_1\\w_1\\q_1},
\mat(1){z_2\\w_2\\q_2}
}_{T_{(y,u,p)}F}:=
\scp{z_1,z_2}_u+G_u(w_1,w_2)+\scp{q_1,q_2}_u
\end{equation}
and thus also $\left(F,\scp{.,.}_{TF}\right)$ as Riemannian manifold.
This scalar product can be used to apply the Riesz representation theorem in order to define the gradient of the Lagrangian
$\grad\LL\in TF$ by the condition
\[
\scp{\grad\LL,h}_{T_{(y,u,p)}F}:= D\LL(y,u,p)h\, , \ \forall h\in  T_{(y,u,p)}F.
\] 
Now, similar to standard nonlinear programming we can solve the problem of finding $(y,u,p)\in F$ with 
\begin{equation}\label{Newton-problem}
\grad\LL(y,u,p)=0
\end{equation}
as a means to find solutions to the optimization problem (\ref{vb-opt1}--\ref{vb-opt2}). The nonlinear problem
(\ref{Newton-problem}) has exactly the form of the root finding problems discussed in \cite{VHS-shape-Riemann}.
Exploiting the Riemannian structure on $TF$, we can formulate a Newton iteration involving the
Riemannian Hessian which is based on the resulting Riemannian connection:
\medskip
\begin{quote}
\fbox{\parbox{\linewidth}{
$k-$th iteration:
\smallskip
\begin{compactenum}[(1)]
\item compute increment $\Delta \xi$ as solution of
\begin{equation}\label{newton}
\hess{}\LL(\xi^k)\Delta \xi=-\grad\LL(\xi^k)
\end{equation}
\item increment $\xi^{k+1}:=\exp_{\xi^k}(\alpha^k\cdot\Delta \xi)$ for some steplength $\alpha^k$
\end{compactenum}
}}
\end{quote}
\medskip
This iteration will be detailed out below. However, before that, we have to specify the scalar product on the Hilbert space involved. Since we will use the exponential map based on the Riemannian metric on $F$, we would like to choose a metric that is in the Hilbert space parts as simple as possible. Therefore we use the metric defined on the Hilbert space $(H_0,\scp{.,.}_0)$ and transfer that canonically to the Hilbert spaces $H(u)$. Thus, we assume now that in the sequel we only have to deal with one particular chart $(U_i,\tau_i)$ from the covering $\{U_i\}$ and define there
\[
\scp{z_1,z_2}_u:=\scp{\tau_i(u)z_1,\tau_i(u)z_2}_0\, , \ \forall u\in U_i\, .
\]
Now, geodesics in the Hilbert space parts of $F$ are represented just by straight lines in $H_0$ and the exponential map can be expressed in the form
\begin{align*}
&\exp_{(y,u,p)}(z,w,q)\\
&=\left(\tau_i\left(\exp^{\cN}_u(w)\right)^{-1}\circ\tau_i(u)(y+z),\exp^{\cN}_u(w),
\tau_i\left(\exp^{\cN}_u(w)\right)^{-1}\circ\tau_i(u)(p+q)\right)
\end{align*}
where $\exp^{\cN}$ denotes the exponential map on the manifold $\cN$.

Within iteration (\ref{newton}), the Hessian has to be discussed. It is based on the 
Riemannian connection $\nabla$ on $F$ at $u\in\cN$. The expression $\nabla^{\cN}_u$ may denote the Riemannian covariant derivative on $T_u\cN$. Since the scalar product in $H$ is completely independent from the location $u\in\cN$, we observe that mixed covariant derivatives of vectors from $H$ with respect to tangential vectors in $T\cN$ are reduced to simple directional derivatives -- which is the case for derivatives in linear spaces anyway. Thus:
\[
\begin{array}{rl}
\nabla_{(h_y,h_u,h_p)}\colon T_{(y,u,p)}F & \hspace{-2mm}\to T_{(y,u,p)}F\\
\mat(1){z\\w\\q} & \hspace{-2mm} \mapsto
\mat(1){
\frac{\partial}{\partial y}z[h_y]+\frac{\partial}{\partial u}z[h_u]+\frac{\partial}{\partial p}z[h_p]\\
\frac{\partial}{\partial y}w[h_y]+\nabla^{\cN}_{u}w[h_u]+\frac{\partial}{\partial p}w[h_p]\\
\frac{\partial}{\partial y}q[h_y]+\frac{\partial}{\partial u}q[h_u]+\frac{\partial}{\partial p}q[h_p]
}
\end{array}
\]
From the definition of the Hessian as $\hess{}\LL[h]:=\nabla_h\grad\LL$ we conlude the following block structure of the Hessian:
\begin{equation}\label{hessian}
\hess{}\LL
=
\mat(3){
D_y\grad_y\LL&D_u\grad_y\LL&D_p\grad_y\LL\\
D_y\grad_u\LL&\nabla^{\cN}_u\grad_u\LL&D_p\grad_u\LL\\
D_y\grad_p\LL&D_u\grad_p\LL&0
}
\end{equation}

\newcommand{\bH}{\bar{H}}
\newcommand{\rH}{\mathring{H}}
\newcommand{\bz}{\bar{z}}
\newcommand{\bw}{\bar{w}}
\newcommand{\bq}{\bar{q}}

From a practical point of view, it may be advantageous to solve equation (\ref{newton}) in a weak formulation as
\begin{equation}\label{weak-form}
\nabla (D\LL(y,u,p)h)\mat(1){z\\w\\q}=-D\LL(y,u,p)h\, ,\quad\forall h\in T_{(y,u,p)}F
\end{equation}
i.e., in detail, the following equations have to be satisfied for all  $h:=(\bz,\bw,\bq)^T\in T_{(y,u,p)}F$:
\begin{alignat}{4}\label{weak-KKT.1}
H_{11}(z,\bz)+H_{12}(w,\bz)+H_{13}(q,\bz)&=-a_u(\bz,p)-\frac{\partial}{\partial y}J(y,u)\bz&\\
\label{weak-KKT.2}
H_{21}(z,\bw)+H_{22}(w,\bw)+H_{23}(q,\bw)&=-\frac{\partial}{\partial u}\left[J(y,u)+a_{u}(y,p)-b_{u}(p)\right]\bw\\
\label{weak-KKT.3}
H_{31}(z,\bq)+H_{32}(w,\bq)+H_{33}(q,\bq)&=-a_{u}(y,\bq)+b_{u}(\bq)
\end{alignat}
where
\begin{alignat*}{2}
H_{11}(z,\bz)&=\frac{\partial^2}{\partial y^2}J(y,u)z\bz\\
H_{12}(w,\bz)&=\frac{\partial}{\partial u}\left[a_u(\bz,p)+\frac{\partial}{\partial y}J(y,u)\bz\right] w\\
H_{13}(q,\bz)&=a_u(\bz,q)\\
H_{21}(z,\bw)&=\frac{\partial}{\partial y}\frac{\partial}{\partial u}\left(\left[J(y,u)+a_{u}(y,p)\right]\bw\right)z\\
H_{22}(w,\bw)&=G\left(\hess{}^{\cN}\left(J(y,u)+a_{u}(y,p)-b_{u}(p)\right)w,\bw\right)\\
H_{23}(q,\bw)&=\frac{\partial}{\partial u}\left[a_{u}(y,q)-b_{u}(q)\right]\bw\\
H_{31}(z,\bq)&=a_{u}(z,\bq)\\
H_{32}(w,\bq)&=\frac{\partial}{\partial u}\left[a_{u}(y,\bq)-b_{u}(\bq)\right]w\\
H_{33}(q,\bq)&=0
\end{alignat*}
One should note that the covariant derivative $\nabla$ reveals natural symmetry properties and thus obvious symmetries can be observed in the components above not involving second shape derivatives. A key observation in \cite{VHS-shape-Riemann} is that even the expression $H_{22}(w,\bw)$ is symmetric in the solution of the shape optimization problem. This motivates a shape--SQP method as outlined below, where away from the solution only expressions in $H_{22}(w,\bw)$ are used which are nonzero at the solution. Its basis is the following observation:

\CROP{
\vspace{2cm}
The algorithm above requires the computation of gradients by the usage of the scalar product in the Hilbert space $H$. Since the bilinear form $a_u$ is typically indefinite or nonsymmetric, it cannot be used as a scalar product in $H$. Thus, in general, additional PDE have to be solved just for the transformation of derivatives into gradients. In order to avoid this, we assume that $y,p\in \rH\subset H\subset\bH$, where $\rH,\bH$ are Hilbert spaces of higher ($\rH$) respectively lower ($\bH$) regularity and that the state constraint and the adjoint equation can be expressed as
\[
G^s(y,u)=0,\qquad G^a(y,u,p)=0
\]
where $G^s:\rH\times\cN$ and $G^a:\rH\times\cN\times\rH\to \bH$ are differentiable. Then, we apply a Newton iteration to the equation
\[
\mat(1){G^a(y,u,p)\\\grad_u\LL(y,u,p)\\G^s(y,u)}=0
\]
where we substitute equation (\ref{newton}) of the algorithm above by
\begin{equation}\label{snewton}
\mat(3){
D_yG^a&\nabla^{\cN}_uG^a&D_pG^a\\
D_y\grad_u\LL&\nabla^{\cN}_u\grad_u\LL&D_p\grad_u\LL\\
D_yG^s&\nabla^{\cN}_uG^s&0
}
\mat(1){z\\w\\q}
=
-
\mat(1){G^a\\\grad_u\LL\\G^s}
\end{equation}
The property of quadratic convergence remains nevertheless, as will be exemplified below. However, the Karush-Kuhn-Tucker matrix is no longer self-adjoint.
}

If the term $H_{22}(w,\bw)$ is replaced by an approximation $\hat{H}_{22}(w,\bw)$, which omits all terms in $H_{22}(w,\bw)$, which are zero at the solution and if the reduced Hessian of (\ref{hessian}) built with this approximation is coercive, equation (\ref{weak-form}) is equivalent to the linear-quadratic problem
\begin{align}\label{qp.1}
&\min\limits_{(z,w)}\frac12\left(H_{11}(z,z)+2H_{12}(w,z)+\hat{H}_{22}(w,w)\right)+\frac{\partial}{\partial y}J(y,u)z +\frac{\partial}{\partial u}J(y,u)w\\\label{qp.2}
&\hspace{0.4cm}\mbox{s.t. }a_u(z,\bar{q})+\frac{\partial}{\partial u}\left[a_{u}(y,\bq)-b_{u}(\bq)\right]w=-a_u(y,\bar{q})+b_u(\bar{q})\, , \ \forall \bar{q}\in H(u)
\end{align}
where the adjoint variable to the constraint (\ref{qp.2}) is just $p+q$. In the next sections, we also omit terms in $H_{11}$ and $H_{12}$, which are zero, when evaluated at the solution of the optimization problems. Nevertheless, quadratic convergence of the resulting SQP method is to be expected and indeed observed in section \ref{sec4}.

\section{Discussion for a Poisson--type model problem}\label{sec3}
In this section, we apply the theoretical discussion of section \ref{sec2} to a PDE constrained shape optimization problem, which is inspired by the standard tracking--type elliptic optimal control problem and motivated by electrical impedance tomography.
It is very close to the model problem of example 2 in \cite{Cea-RAIRO} and the inverse interface problem in \cite{ItoKunisch}. 

\begin{figure}[h]
\label{fig_omega}
\begin{center}
   \begin{overpic}[width=0.3\textwidth]{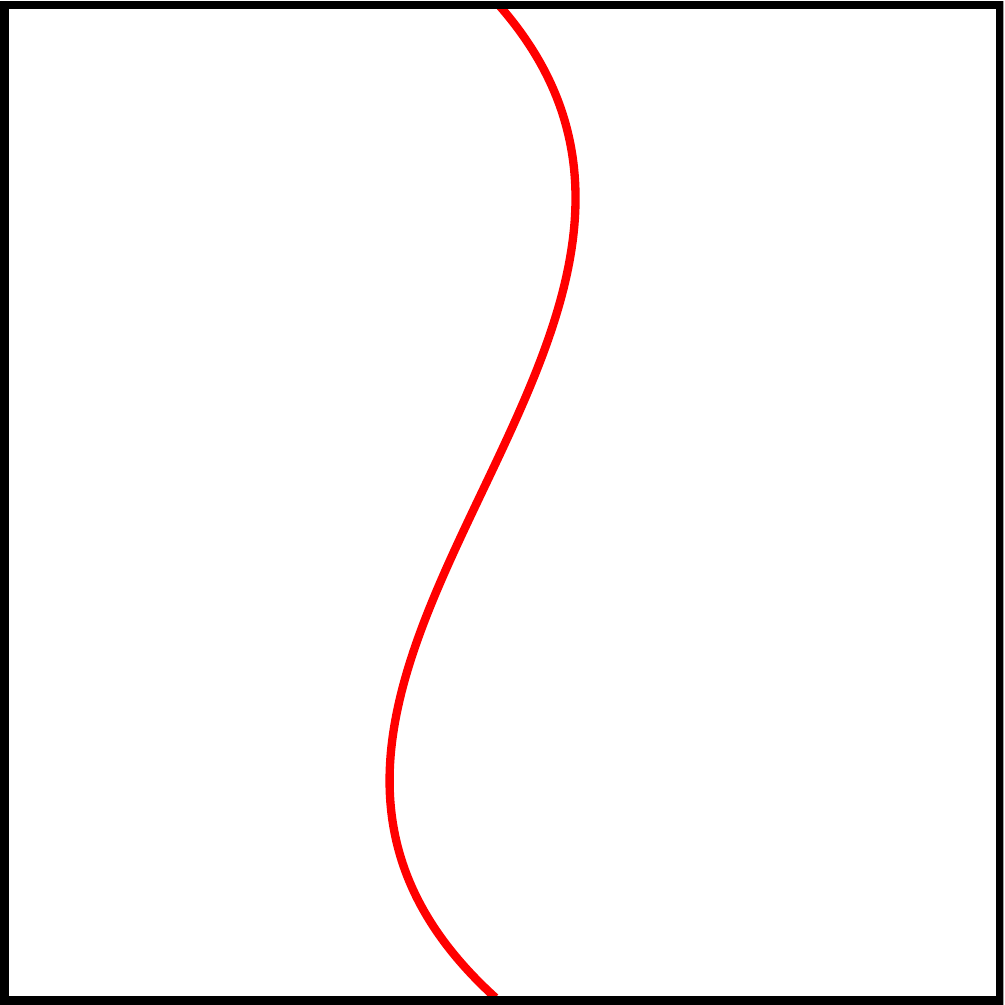}
   \put(70,80){$\Omega_2(u) $}
   \put(10,80){$\Omega_1(u)$}
   \put(43,13){$\color{red}u$}
   \end{overpic}
   \end{center}
   \caption{Example of a domain $\Omega(u)=\Omega_1(u)\cupdot u \cupdot\Omega_2(u)$}
\end{figure}

Let the domain $\Omega:=(0,1)^2\subset \R^2$ split into the two subdomains $\Omega_1,\Omega_2 \subset \Omega$ such that $\Omega_1 \cupdot\Gamma\cupdot\Omega_2 = \Omega$ and $\partial\Omega_1 \cap \partial\Omega_2 = \Gamma$.  The interface $\Gamma$ is replaced by $u$ and an element of the following manifold
\[
B_e^0([0,1],\R^2):=\mbox{Emb}^0([0,1],\R^2)/\mbox{Diff}^0([0,1])
\] 
i.e., an element of the set of all equivalence classes of the set of embeddings
\begin{align*}
\mbox{Emb}^0([0,1],\R^2):= \{\phi\in C^\infty([0,1],\R^2)\, |\, &
\phi(0)=(0.5,0),\phi(1)=(0.5,1),\\
 &\phi \mbox{ injective immersion}\}
\end{align*}
where the equivalence relation is defined by the set of all $C^\infty$ re--parameterizations, i.e., by the set of all diffeomorphisms
\begin{align*}
\mbox{Diff}^0([0,1],\R^2):= \{\phi\colon [0,1]\to[0,1]\, |\, &
\phi(0)=(0.5,0),\phi(1)=(0.5,1),\\
 &\phi \mbox{ diffeomorphism}\}.
\end{align*}
In Figure \ref{fig_omega} the construction of the domain $\Omega$ from the interface $u\in B_e^0([0,1],\R^2)$ is illustrated. Now, we consider $\Omega$ dependent on $u$. Therefore, we denote it by $\Omega(u)=\Omega_1(u)\cupdot u \cupdot\Omega_2(u)$. 

\begin{remark}
The manifold  $B_e^0([0,1],\R^2)$ is constructed in analogy to the manifold $B_e(S^1,\R^2)$ in \cite{MM-2006} as a set of equivalence classes in a set of embeddings with respect to a equivalence relation which is given by a set of diffeomorphisms.
Moreover, a particular point on the manifold $B_e^0([0,1],\R^2)$ is represented by a curve $c\colon [0,1]\to\R^2\text{, }\theta\mapsto c(\theta)$. Because of the equivalence relation $\mbox{Diff}([0,1])$, the tangent space is isomorphic to the set of all smooth vector fields along $c$, i.e.,
$$T_cB_e^0([0,1],\R^2)\cong\{h\, |\,h=\alpha n, \alpha\in C^\infty \left([0,1],\R\right)\}$$
where n is the unit outer normal to $\Omega_1(u)$ at $u$. Thus, all considerations of \cite{VHS-shape-Riemann} carry easily over to our manifold $B_e^0([0,1],\R^2)$.
\end{remark}
 
The PDE constrained shape optimization problem is given in strong form by
\begin{align}\label{oc1}
&\min_{u} \hspace{0,1cm} J(y,u)\equiv\frac12\int_{\Omega(u)}(y-\bar{y})^2dx+
\mu\int_u1ds\\
\label{oc2}
&\hspace{0,2cm}\mbox{s.t. }- \triangle y=f \quad\text{in } \Omega(u)
\\
\label{oc3}
&\hspace{15.7mm}y=0\quad\text{on }\partial\Omega(u)
\end{align}
where
\begin{equation}
f\equiv\begin{cases}
f_1 = \text{const.}\quad\text{in }\Omega_1(u)\\
f_2 = \text{const.} \quad\text{in }\Omega_2(u)
\end{cases}
\hspace{-.2cm}.
\end{equation}
The perimeter regularization with $\mu>0$ in the objective (\ref{oc1}) is a frequently used means to overcome ill--posedness of the optimization problem  (e.g.~\cite{Burger-2004}).
Let $n$ be the unit outer normal to $\Omega_1(u)$ at $u$. We observe that the unit outer normal to $\Omega_2(u)$ at $u$ is equal to $-n$, which enables us to use only one normal n for the subsequent discussions. Furthermore, we have interface conditions at the interface $u$. We formulate explicitly the continuity of the state and of the flux at the boundary $u$ as
\begin{align}
\label{n}
\left\llbracket y \right\rrbracket =0\text{, }\quad\left\llbracket\frac{\partial y}{\partial n}\right\rrbracket=0\quad\text{on }u
\end{align}
where the jump symbol $\left\llbracket\cdot\right\rrbracket$ denotes the discontinuity across the interface $u$ and is defined by $\left\llbracket v\right\rrbracket:=v_1-v_2$ where $v_1:=v\,\rule[-2mm]{.1mm}{4mm}_{\hspace{.6mm}\Omega_1}$ and $v_2:=v\,\rule[-2mm]{.1mm}{4mm}_{\hspace{.6mm}\Omega_2}$.

The boundary value problem (\ref{oc2}-\ref{n}) is written in weak form as
\begin{equation}
\label{wf}
a_{u}(y,p)=b_{u}(p)\, , \ 
\forall p\in H_0^1(\Omega(u))
\end{equation}
where
\begin{align}\label{sweak.1}
a_{u}(y,p)&:=
\int_{\Omega(u)}\nabla y^T\nabla p\hspace{.5mm}dx-\int_u \left\llbracket \frac{\partial y}{\partial n}p\right\rrbracket ds\\
\label{sweak.2}
b_{u}(p)&:=\int_{\Omega(u)}fp\hspace{.5mm}dx.
\end{align}

Now, $F$ from definition \ref{lagrangian} takes the specific form
\[
F:=\left\{(H_0^1(\Omega(u)),u,H_0^1(\Omega(u)))\,|\,u\in B_e^0([0,1],\R^2)\right\}.
\]
The metric in the vector space parts is constructed by employing a "mesh deformation". Mesh deformations are often used to deform a computational mesh smoothly in accordance with a deformation of the boundary of the computational domain. Here, we use this in the form of a deformation of the computational domain rather than of the mesh only and assume that there is a bijective $C^\infty$--mapping
\[
\Phi_u\colon[0,1]^2\to \Omega(u),
\]
e.g., $\Phi_u$ is the deformation given by the solution of a linear elasticity problem.
Thus, we can construct the necessary bijective identification
\begin{equation*}
\tau(u)\colon H_0^1(\Omega(u))\to H_0^1\left((0,1)^2\right)\text{, }g\mapsto g\circ\Phi_u.
\end{equation*} 

We have to detail the expressions in equation (\ref{newton}) or respectively (\ref{weak-form}). First, the Lagrangian is defined for $(y,u,p)\in F$ as
\[
\LL(y,u,p):=J(y,u)+
a_u(y,p)-b_u(p)
\]
where $J(y,u)$ is defined in (\ref{oc1}) and $a_u, b_u$ are defined in (\ref{sweak.1}, \ref{sweak.2}). Now, we focus on the shape derivative of $\LL$ in direction of a continuous vector field $V$. It is defined by
\begin{equation}
\label{def_shape_der}
\frac{\partial}{\partial u}\LL(y,u,p)[V]:= \lim\limits_{t\to 0^+}\frac{\LL(y,u_t,p)-\LL(y,u,p)}{t}
\end{equation}
if for all directions $V$ this limit exists and the mapping $V\mapsto \frac{\partial}{\partial u}\LL(y,u,p)[V]$ is linear and continuous. The perturbed boundaries $u_t$ in (\ref{def_shape_der}) are defined by
\begin{equation}
u_t:=F_t(u)=\{F_t(x)\colon x\in u\}\text{ with }u_0=u
\end{equation}
where $F_t(x):=x+tV(x)$ denotes the perturbation of identity and $t\in[0,T]$ with $T>0$. 
\begin{remark}
One should note that we get perturbed domains $\Omega_t$ given by
\begin{equation}
\Omega_t:=F_t(\Omega(u))=\{F_t(x)\colon x\in \Omega(u)\}\text{ with }\Omega_0=\Omega(u)
\end{equation}
due to the perturbed boundaries $u_t$.
\end{remark}
\begin{remark}
The perturbation of $u$ or respectively $\Omega(u)$ could also be described by the velocity method, i.e., as the flow $F_t(x):=\xi(t,x)$ determined by the initial value problem
\begin{equation}
\begin{split}
\frac{d\xi(t,x)}{dt}&=V(\xi(t,x))\\
\xi(0,x)&=x
\end{split}
\end{equation}
instead of the perturbation of identity.
\end{remark}

We first consider the objective $J$ in (\ref{oc1}) without perimeter regularization. Then the shape derivative $\frac{\partial}{\partial u}\LL(y,u,p)[V]$ can be expressed as an integral over the domain $\Omega(u)$ as well as an integral over the interface $u$ which is better suited for a finite element implementation as already mentioned for example in \cite[remark 2.3, p. 531]{Delfour-Zolesio-2001}. An important point to note here is that the shape derivative of our $\LL$ evaluated in its saddle--point is equal to the one of $J$ due to the theorem of Correa and Seeger \cite[theorem 2.1]{CorreaSeeger}. Such a saddle--point is given by
\begin{align}
\frac{\partial \LL(\Omega,y,p)}{\partial y}=\frac{\partial \LL(\Omega,y,p)}{\partial p}=0\label{saddlepointcond}
\end{align}
which leads to the adjoint equation
\begin{align}
- \triangle p&=-(y-\overline{y}) \quad\text{in } \Omega(u)\label{adjoint1}\\
p&=0\quad\text{on }\partial\Omega(u)\\
\left\llbracket p \right\rrbracket &=0\quad\text{on } u\label{adjoint4}\\
\left\llbracket\frac{\partial p}{\partial n}\right\rrbracket&=0\quad\text{on }u\label{adjoint3}
\end{align}
and to the state equation
\begin{equation}
- \triangle y=f \quad\text{in } \Omega(u).\label{designe}
\end{equation}
Like in \cite{Interface_parabolic} we first deduce a representation of the shape derivative expressed as a domain integral which will later allow us to calculate the boundary expression of the shape derivative by means of integration by parts on the interface $u$. One should note however, that by the Hadamard structure theorem \cite[theorem 2.27]{Sokolowski-1992} only the normal part of the continuous vector field has an impact on its value. Applying the following common rule for differentiating domain integrals
\begin{equation}
\label{der_domain_int}
\frac{d^+}{dt}\left(\int_{\Omega_t}\eta(t)\right)\,\rule[-4mm]{.1mm}{9mm}_{\hspace{1mm}t=0}=\int_{\Omega}\left(D_m\eta+\mathrm{div}(V)\eta\right)
\end{equation}
which was proved in \cite[lemma 3.3]{HaslingerMakinen} yields
\begin{equation}
\label{shape_der_1}
\begin{split}
\frac{\partial}{\partial u}\LL(y,u,p)[V]
&=\lim\limits_{t\to 0^+}\frac{\LL(y,u_t,p)-\LL(y,u,p)}{t}=\frac{d^+}{dt}\LL(y,u_t,p)\,\rule[-3mm]{.1mm}{6mm}_{\hspace{.5mm}t=0}\\
&=\int_{\Omega(u)}D_m\left(\frac{1}{2}(y-\overline{y})^2\right)+D_m\left(\nabla y^T\nabla p\right)-D_m(fp)\\
&\hspace*{13mm}+\mathrm{div}(V)\left(\frac{1}{2}(y-\overline{y})^2+\nabla y^T\nabla p-fp\right)\hspace{.7mm}dx\\
&\hspace*{4mm}-\int_u D_m\left( \left\llbracket \frac{\partial y}{\partial n} p\right\rrbracket\right) + \mathrm{div}_u(V)\left\llbracket \frac{\partial y}{\partial n} p\right\rrbracket\hspace{.7mm}ds
\end{split}
\end{equation}
where $D_m$ denotes the material derivative with respect to $F_t=id+tV$ which is defined by
\begin{equation}
\label{material}
D_m\left(j(x)\right):=\lim\limits_{t\to 0^+}\frac{\left(j\circ F_t\right)(x)-j(x)}{t}=\frac{d^+}{dt}\left(j\circ F_t\right)(x)\,\rule[-2.5mm]{.1mm}{6mm}_{\hspace{1mm}t=0}
\end{equation}
for a generic function $j\colon \Omega_t\to\R$. For the material derivative the product rule holds. Moreover, the following equality was proved in \cite{Berggren}
\begin{equation}
\label{material_grad}
D_m\left(\nabla j\right)=\nabla \left(D_m(j)\right)-\nabla V^T\nabla j.
\end{equation}
Combining (\ref{shape_der_1}), the product rule and (\ref{material_grad}) we obtain
\begin{equation}
\label{shape_der_2}
\begin{split}
\frac{\partial}{\partial u}\LL(y,u,p)[V]
&=\int_{\Omega(u)}(y-\overline{y})D_m\left(y\right)+\nabla\hspace{-.5mm} \left(D_m(y)\right)^T\nabla p+\nabla y^T\nabla\hspace{-.5mm} \left(D_m\left(p\right)\right)\\
&\hspace*{13mm}-\nabla y^T \left(\nabla V +\nabla V^T\right)\nabla p-D_m(f)p-fD_m(p)\\
&\hspace*{13mm}+\mathrm{div}(V)\left(\frac{1}{2}(y-\overline{y})^2+\nabla y^T\nabla p-fp\right)\hspace{.7mm}dx\\
&\hspace*{4mm}-\int_u  \left\llbracket D_m\left(\frac{\partial y}{\partial n}\right) p+\frac{\partial y}{\partial n}D_m(p)\right\rrbracket + \mathrm{div}_u(V)\left\llbracket \frac{\partial y}{\partial n} p\right\rrbracket\hspace{.7mm}ds.
\end{split}
\end{equation}
From this we get
\begin{equation}
\label{shape_der_3}
\begin{split}
\frac{\partial}{\partial u}\LL(y,u,p)[V]
&=\int_{\Omega(u)}\left((y-\overline{y})-\triangle p\right)D_m\left(y\right)+\left(-\triangle y-f\right)D_m\left(p\right)\\
&\hspace*{13mm}-\nabla y^T \left(\nabla V +\nabla V^T\right)\nabla p-D_m(f)p\\
&\hspace*{13mm}+\mathrm{div}(V)\left(\frac{1}{2}(y-\overline{y})^2+\nabla y^T\nabla p-fp\right)\hspace{.7mm}dx\\
&\hspace*{4mm}+\int_u  \left\llbracket\frac{\partial p}{\partial n}D_m(y) -D_m\left(\frac{\partial y}{\partial n}\right) p\right\rrbracket + \mathrm{div}_u(V)\left\llbracket \frac{\partial y}{\partial n} p\right\rrbracket\hspace{.7mm}ds.
\end{split}
\end{equation}
To deal with the term $D_m(f)p$, we note that the shape derivative of a generic function $j\colon \Omega_t\to\R$ with respect to the vector field $V$ is given by
\begin{equation}
\label{shape_material_der}
Dj[V]:=D_mj-V^Tj.
\end{equation}
Therefore $D_m(f)p$ is equal to $pV^T\nabla f$ in the both subdomains $\Omega_1(u)$, $\Omega_2(u)$. Due to the continuity of the state and of the flux (\ref{n}) their material derivative is continuous. Thus, we get
\begin{align}
\left\llbracket \frac{\partial p}{\partial n}D_m(y)\right\rrbracket&=D_m(y)\left\llbracket \frac{\partial p}{\partial n}\right\rrbracket\stackrel{(\ref{adjoint3})}=0\quad \text{on }u\\
\left\llbracket D_m\hspace{-.5mm}\left(\frac{\partial y}{\partial n}\right)p\right\rrbracket&=D_m\hspace{-.5mm}\left(\frac{\partial y}{\partial n}\right)\left\llbracket p\right\rrbracket\stackrel{(\ref{adjoint4})}=0\quad \text{on }u.
\end{align}
That
\begin{equation}
\label{pbe}
\left\llbracket \frac{\partial y}{\partial n}p\right\rrbracket=0\quad \text{on }u
\end{equation}
follows from (\ref{n}), (\ref{adjoint4}) and the identity
\begin{equation}
\left\llbracket ab\right\rrbracket=\left\llbracket a\right\rrbracket b_1 +a_2 \left\llbracket b\right\rrbracket= a_1 \left\llbracket b\right\rrbracket+\left\llbracket a\right\rrbracket b_2
\end{equation}
which implies
\begin{equation}
\label{bracket0}
\left\llbracket ab\right\rrbracket=0 \text{ if } \left\llbracket a\right\rrbracket=0\wedge \left\llbracket b\right\rrbracket=0.
\end{equation}
By combining (\ref{adjoint1}), (\ref{designe}) and (\ref{shape_der_3}--\ref{pbe}), we obtain 
\begin{equation}
\label{boundary_expression}
\boxed{
\begin{split}
\frac{\partial}{\partial u}\LL(y,u,p)[V]=\int_{\Omega(u)}&-\nabla y^T \left(\nabla V +\nabla V^T\right)\nabla p-pV^T\nabla f\\
&+\mathrm{div}(V)\left(\frac{1}{2}(y-\overline{y})^2+\nabla y^T\nabla p-fp\right)\hspace{.7mm}dx
\end{split}
}
\end{equation}
i.e., the shape derivative of $\LL$ expressed as domain integral which is equal to the one of $J$ due to the theorem of Correa and Seeger. Now, we convert this domain integral into a boundary integral as mentioned above. Integration by parts in (\ref{boundary_expression}) yields
\begin{equation}
\begin{split}
&\int_{\Omega(u)} \mathrm{div}(V)\left(\frac{1}{2}(y-\overline{y})^2+\nabla y^T\nabla p-fp\right)dx\\
&=-\int_{\Omega(u)}V^T\left((y-\overline{y})\nabla y+\nabla\left(\nabla y^T\nabla p\right)-\nabla (f p)\right)\hspace{.3mm}dx\\
&\hspace*{3.5mm}+\int_{u} \left\llbracket\left(\frac{1}{2}(y-\overline{y})^2+\nabla y^T\nabla p-fp\right)\left<V,n\right>\right\rrbracket ds\\
&\hspace*{3.5mm}+\int_{\partial\Omega(u)}\left(\frac{1}{2}(y-\overline{y})^2+\nabla y^T\nabla p-fp\right)\left<V,n\right>ds.\label{int_by_parts1}
\end{split}
\end{equation}
Since the outer boundary $\partial\Omega$ is not variable, we can choose the deformation vector field $V$ equals zero in small neighbourhoods of $\partial\Omega(u)$. Therefore, the outer integral in (\ref{int_by_parts1}) disappears. Combining (\ref{boundary_expression}), (\ref{int_by_parts1}) and the vector calculus identity 
\begin{equation*}
\nabla y^T\left(\nabla V+\nabla V^T\right)\nabla p+V^T\nabla\left(\nabla y^T\nabla p\right)=\nabla p^T\nabla\left(V^T\nabla y\right)+\nabla y^T \nabla\left(V^T\nabla p\right)
\end{equation*}
which was proved in \cite{Berggren} gives
\begin{equation}
\begin{split}
\frac{\partial}{\partial u}\LL(y,u,p)[V]=&\int_{\Omega(u)}-\nabla p^T\nabla\left(V^T\nabla y\right)-\nabla y^T \nabla\left(V^T\nabla p\right)\\
&\hspace*{9mm}-(y-\overline{y})V^T\nabla y+fV^T\nabla p \hspace{.7mm}dx\\
&+\int_{u} \left\llbracket\left(\frac{1}{2}(y-\overline{y})^2+\nabla y^T\nabla p-fp\right)\left<V,n\right>\right\rrbracket ds. \label{sg1}
\end{split}
\end{equation}
Then, applying integration by parts in (\ref{sg1}) we get
\begin{equation}
\begin{split}
&\int_{\Omega(u)}\nabla y^T\nabla\left(V^T\nabla p\right)dx\\
&=-\int_{\Omega(u)}\triangle yV^T\nabla p\hspace{.7mm}dx+\int_{u}\left\llbracket \frac{\partial y}{\partial n}V^T\nabla p \hspace{.3mm}\right\rrbracket ds+\int_{\partial\Omega(u)}\frac{\partial y}{\partial n}V^T\nabla p \hspace{.7mm}ds\label{int_by_parts2}
\end{split}
\end{equation}
and analogously
\begin{equation}
\begin{split}
& \int_{\Omega(u)}\nabla p^T\nabla\left(V^T\nabla y\right)dx\\
&=-\int_{\Omega(u)}\triangle pV^T\nabla y\hspace{.7mm}dx+\int_{u}\left\llbracket \frac{\partial p}{\partial n}V^T\nabla y \hspace{.3mm}\right\rrbracket ds+\int_{\partial\Omega(u)}\frac{\partial p}{\partial n}V^T\nabla y \hspace{.7mm}ds. \label{int_by_parts3}
\end{split}
\end{equation}
Like in (\ref{int_by_parts1}) the outer integral in (\ref{int_by_parts2}) as well as in (\ref{int_by_parts3}) vanishes due to the fixed outer boundary $\partial\Omega(u)$. Thus, it follows that
\begin{equation}
\begin{split}
\frac{\partial}{\partial u}\LL(y,u,p)[V]=&\int_{\Omega(u)}V^T\nabla p\left(\triangle y+f\right)+V^T\nabla y\left(\triangle p -(y-\overline{y})\right)\hspace{.7mm}dx\\
&+\int_u \left\llbracket\left(\frac{1}{2}(y-\overline{y})^2+\nabla y^T\nabla p-fp\right)\left<V,n\right>\right\rrbracket\\
&\hspace*{9mm}-\left\llbracket \frac{\partial y}{\partial n}V^T\nabla p \hspace{.3mm}\right\rrbracket-\left\llbracket \frac{\partial p}{\partial n}V^T\nabla y \hspace{.3mm}\right\rrbracket\hspace{.7mm}ds.\label{sg2}
\end{split}
\end{equation}
The domain integral in (\ref{sg2}) vanishes due to (\ref{adjoint1}) and (\ref{designe}). Moreover, the term $\left\llbracket \frac{1}{2}(y-\overline{y})^2\left<V,n\right>\right\rrbracket$ disappears because of (\ref{n}) and the term $\left\llbracket \nabla y^T\nabla p\left<V,n\right>\right\rrbracket$ because of the continuity of $\nabla y$ and $\nabla p$. That
\begin{equation}
\label{pb}
\left\llbracket \frac{\partial y}{\partial n}V^T\nabla p\right\rrbracket=\left\llbracket \frac{\partial p}{\partial n}V^T\nabla y\right\rrbracket =\left<V,n\right>\left\llbracket \frac{\partial y}{\partial n}\frac{\partial p}{\partial n}\right\rrbracket=0
\end{equation}
follows from (\ref{n}), (\ref{adjoint3}) and (\ref{bracket0}). Thus, we obtain the shape derivative of $\LL$ expressed as interface integral:
\begin{equation}
\label{shape_der_interface}
\boxed{
\frac{\partial}{\partial u}\LL(y,u,p)[V]=-\int_u \left\llbracket f\right\rrbracket p\left<V,n\right>ds
}
\end{equation}

Now, we consider the objective $J$ in (\ref{oc1}) with perimeter regularization. Combining (\ref{shape_der_interface}) with proposition 5.1 in \cite{Novruzi-2002} we get
\begin{equation}
\label{shape_derivative_interface_reg}
\boxed{
\frac{\partial}{\partial u}\LL(y,u,p)[V]=\int_u \left(-\left\llbracket f\right\rrbracket p+\mu\kappa\right)\left<V,n\right>ds
}
\end{equation}
where $\kappa$ denotes the curvature corresponding to the normal $n$.

\begin{remark}
Note that (\ref{shape_der_interface}) is equal to $\frac{\partial}{\partial u}J(y,u)[V]$ without perimeter regularization and (\ref{shape_derivative_interface_reg}) is equal to $\frac{\partial}{\partial u}J(y,u)[V]$ with perimeter regularization due to the theorem of Correa and Seeger as mentioned above.
\end{remark}

We focus now on the weak formulation (\ref{weak-KKT.1}-\ref{weak-KKT.3}) and observe the following for the right hand sides in the case of (\ref{oc1}--\ref{oc3}):
\begin{align}
-a_u(\bz,p)-\frac{\partial}{\partial y}J(y,u)\bz =& -\int_{\Omega(u)}\nabla \bz^T\nabla p+(y-\bar{y})\bz \hspace{.7mm}dx\label{Diri.adj}\\
-\frac{\partial}{\partial u}\left[J(y,u)+a_{u}(y,p)-b_{u}(p)\right]\bw =&\int_u \left(\left\llbracket f\right\rrbracket p-\mu\kappa\right)\left<\bw,n\right>ds\label{no2}\\
-a_{u}(y,\bq)+b_{u}(\bq) =&\int_{\Omega(u)}-\nabla y^T\nabla \bq +f\bq \hspace{.7mm}dx\label{Diri.state}
\end{align}
These expressions are set to zero, in order to define the necessary conditions of optimality. 

Now, we discuss more details about the Hessian operators in the left hand sides of (\ref{weak-KKT.1}--\ref{weak-KKT.3}). We first consider them without the term $H_{22}$ which requires special care. These are at the solution $(y,u,p) \in F$ of the optimization problem (\ref{oc1}--\ref{oc3}) for all $h:=(\bz,\bw,\bq)^T\in T_{(y,u,p)}F$ as follows:
\begin{alignat*}{2}
H_{11}(z,\bz)&= \frac{\partial^2}{\partial y^2}J(y,u)z\bz=\int_{\Omega(u)}z\bz dx\\
H_{12}(w,\bz)&=\frac{\partial}{\partial u}\left[a_u(\bz,p)+\frac{\partial}{\partial y}J(y,u)\bz\right] w=0\\
H_{13}(q,\bz)&=a_u(\bz,q)=\int_{\Omega(u)}\nabla \bz^T \nabla q dx\\
H_{21}(z,\bw)&=\frac{\partial}{\partial y}\frac{\partial}{\partial u}(\left[J(y,u)+a_{u}(y,p)\right]\bw)z=0\\
H_{23}(q,\bw)&=\frac{\partial}{\partial u}\left[a_{u}(y,q)-b_{u}(q)\right]\bw=-\int_u\left\llbracket f\right\rrbracket q\left<\bw,n\right>ds\\
H_{31}(z,\bq)&=a_{u}(z,\bq)=\int_{\Omega(u)}\nabla z^T \nabla \bq dx\\
H_{32}(w,\bq)&=\frac{\partial}{\partial u}\left[a_{u}(y,\bq)-b_{u}(\bq)\right]w=-\int_u\left\llbracket f\right\rrbracket \bq\left<w,n\right>ds\\
H_{33}(q,\bq)&=0
\end{alignat*}
We compute now the term $H_{22}$. It will be evaluated at the solution of the optimization problem which means that it consists only of the second shape derivative. In section \ref{sec4} this solution will be a straight line connection of the points $(0.5,0)$ and $(0.5,1)$, i.e., the curvature is equal to zero. 
Combining proposition 5.1 in \cite{Novruzi-2002} with the following rule for differentiating boundary integrals
\begin{equation}
\label{der_bound_int}
\frac{d^+}{dt}\left(\int_{\Gamma_t}\eta(t)\right)\,\rule[-4mm]{.1mm}{9mm}_{\hspace{1mm}t=0}=\int_{\Gamma} \left(D\eta[V]+\left(\frac{\partial\eta}{\partial n}+\eta\kappa\right)\left<V,n\right>\right)
\end{equation}
which was proved in \cite{HCK86} yields
\begin{equation}
\label{H22_1}
\begin{split}
& H_{22}(w,\bw)\\
&=G(\hess{}^{\cN}\left(J(y,u)+a_{u}(y,p)-b_{u}(p)\right)w,\bw)\\
&= \int_u -D\left(\left\llbracket f\right\rrbracket p\right)[\bw]\left<w,n\right>-\left\llbracket f\right\rrbracket\left(\kappa p+\frac{\partial p}{\partial n}\right)\left<\bw,n\right>\left<w,n\right>\\
&\hspace*{8mm}+\mu\frac{\partial w}{\partial \tau}\frac{\partial \bar{w}}{\partial \tau}\left<\bw,n\right>\left<w,n\right>ds
\end{split}
\end{equation}
where $\partial/\partial \tau$ denotes the derivative tangential to $u$. We have to evaluate the shape derivative $D\left(\left\llbracket f\right\rrbracket p\right)[\bw]$ in (\ref{H22_1}). We observe in our special case
\begin{equation}
\label{p0}
p=0\quad\text{on }u
\end{equation}
because of the necessary optimality condition (\ref{no2}). Thus, it holds that
\begin{equation}
\label{Dp}
Dp[\bw]=-\bw^T\nabla p=-\bw^T\frac{\partial p}{\partial n}n\quad\text{on }u
\end{equation}
due to (\ref{shape_material_der}). Applying the product rule for shape derivatives yields
\begin{equation}
\label{Dfp}
\begin{split}
D\left(\left\llbracket f\right\rrbracket p\right)[\bw] & = \left\llbracket Df[\bw]\hspace{.7mm}p\right\rrbracket +\left\llbracket f Dp[\bw]\right\rrbracket\stackrel{(\ref{n})}=\left\llbracket Df[\bw]\right\rrbracket p+\left\llbracket f \right\rrbracket Dp[\bw]\\
&\hspace*{-1.8mm}\underset{(\ref{Dp})}{\stackrel{(\ref{p0})}=}-\left\llbracket f \right\rrbracket \frac{\partial p}{\partial n}\left<\bw,n\right>\quad\text{on }u.
\end{split}
\end{equation}
Thus, the Hessian operator $H_{22}$ reduces to
\begin{equation}
\hat{H}_{22}(w,\bw) =\int_u \left(\mu\frac{\partial w}{\partial \tau}\frac{\partial \bar{w}}{\partial \tau}-\left\llbracket f\right\rrbracket\kappa p\right)\left<w,n\right>\left<\bw,n\right>ds.
\end{equation}

%

By using the expressions above, we can formulate the QP (\ref{qp.1}, \ref{qp.2}) at the solution in the following form:

\begin{align}\label{qp.p1}
&\hspace*{-1.5mm}\min\limits_{(z,w)}F(z,w,y,p) \\
&\mbox{s.t. }
\int_{\Omega(u)}\nabla z^T\nabla \bq \hspace{.7mm}dx-\int_u\left\llbracket f\right\rrbracket\bq w \hspace{.7mm}ds \nonumber \\\label{qp.p2}
& \hspace*{1.2cm}=-\int_{\Omega(u)}\nabla y^T\nabla \bq \hspace{.7mm}dx +\int_{\Omega(u)}f\bq \hspace{.7mm}dx\, , \ \forall \bar{q}\in H_0^1(\Omega(u))
\end{align}
where the objective function $F$ is given by
\begin{align}
F(z,w,y,p)=&\int_{\Omega(u)}\frac{z^2}2+ (y -\bar{y}) z\hspace{.7mm}dx+ \int_u \mu\kappa w -\left\llbracket f\right\rrbracket pw\hspace{.7mm}ds\nonumber\\
&+\frac{1}{2} \int_u \mu\left(\frac{\partial w}{\partial\tau}\right)^2-\left\llbracket f\right\rrbracket \kappa p w^2\hspace{.7mm}ds.
\end{align}
This QP in weak formulation can be rewritten in the more intelligible strong form of an optimal control problem:

\begin{align}
&\label{qps.1}
\hspace*{-17.5mm}\min\limits_{(z,w)}F(z,w,y,p)\\
\label{qps.2}\mbox{s.t. } -\triangle z &= \triangle y +f_1 \quad \mbox{in }\Omega_1(u)\\
\label{qps.3}-\triangle z &= \triangle y +f_2\quad\mbox{in }\Omega_2(u)\\
\label{qps.4}\frac{\partial z}{\partial n}&=f_1 w \quad\mbox{on } u\\
\label{qps.5}-\frac{\partial z}{\partial n}&=f_2 w\quad\mbox{on } u\\
\label{qps.6}z&=0\quad \mbox{on }\partial\Omega(u)
\end{align}

The adjoint problem to this optimal control problem is the boundary value problem:
\begin{align}\label{ocad.1}
-\triangle q &= -z-(y-\bar{y}) \quad \mbox{in }\Omega(u)\\
\label{ocad.2}
           q &= 0 \quad \mbox{on }\partial\Omega(u)
\end{align}
The resulting design equation for the optimal control problem (\ref{qps.1}--\ref{qps.6}) is 
\begin{equation}
\label{oc.design}
0= -\left\llbracket f\right\rrbracket\left(p+\kappa pw+ q\right)+\mu\kappa-\mu\frac{\partial^2 w}{\partial\tau^2}
\quad \mbox{on }u.
\end{equation}

\section{Numerical Results}\label{sec4}
\begin{figure}[h]
\label{fig-iterations}
\unitlength1cm
\begin{picture}(12,4)
\put(0,0){\includegraphics[scale=0.05]{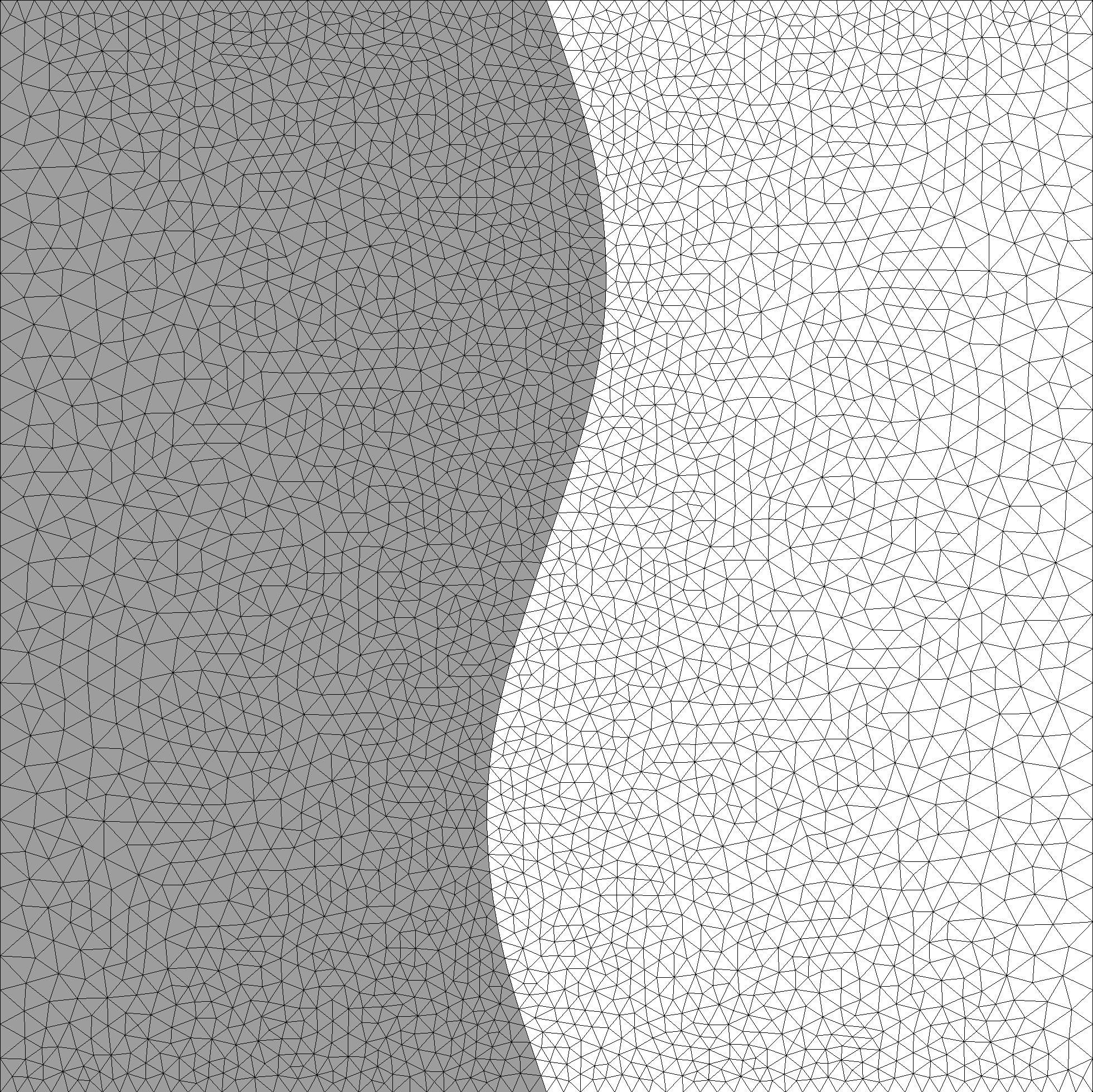}}
\put(4.4,0){\includegraphics[scale=0.05]{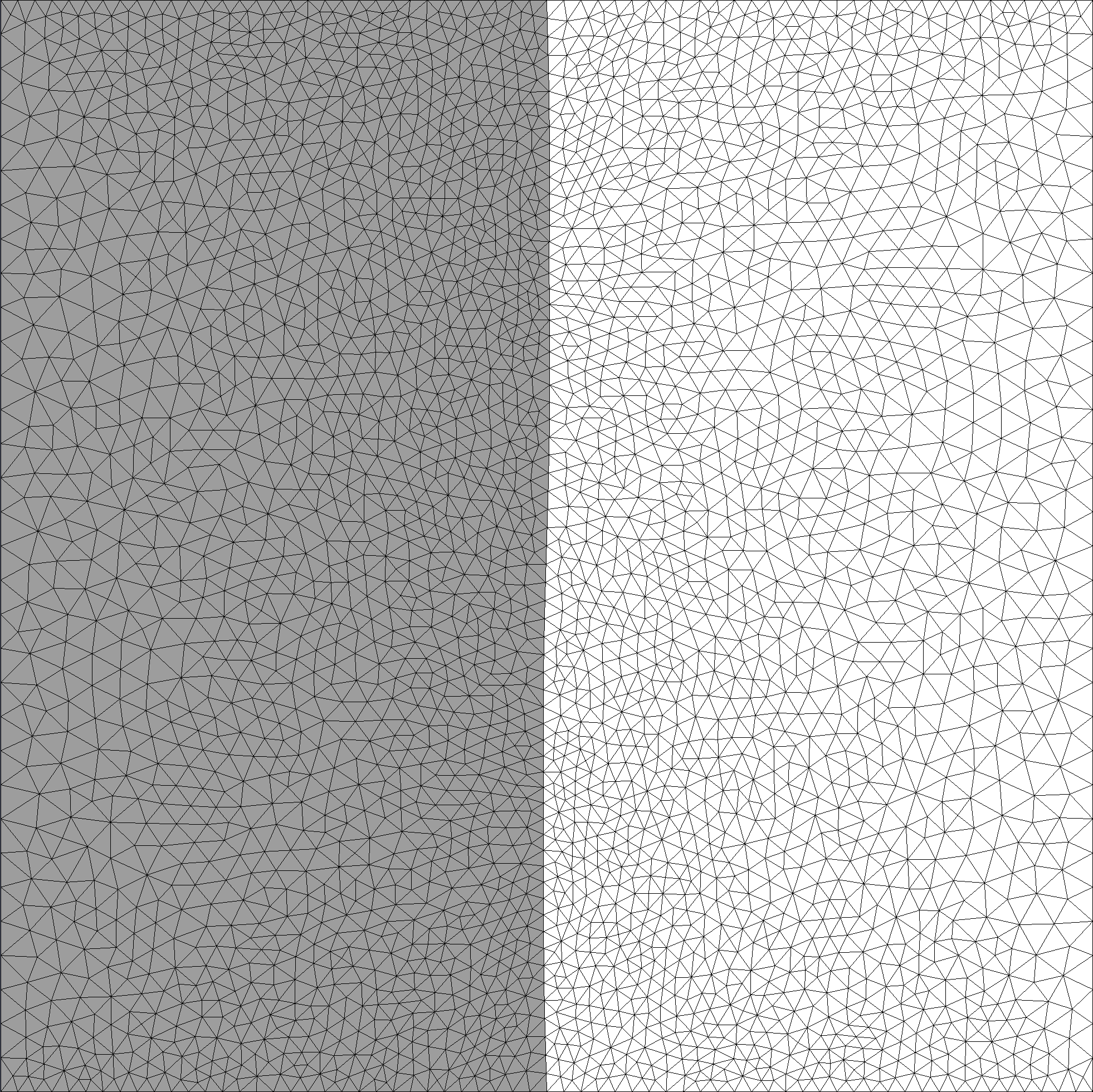}}
\put(8.8,0){\includegraphics[scale=0.05]{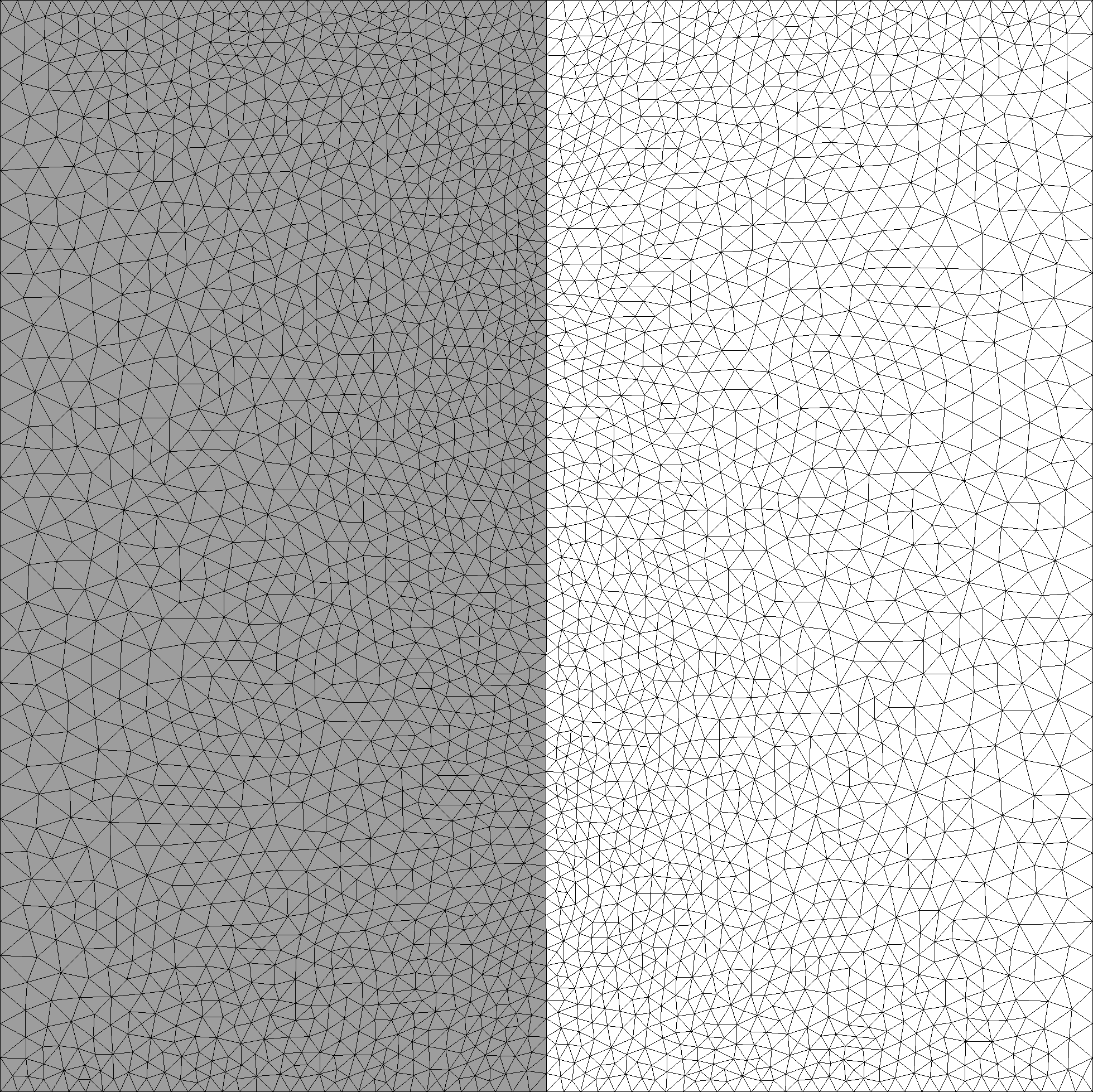}}
\end{picture}
\caption{Iterations 0, 1 and 2 (left to right) together with deformations of mesh $\Omega_h^1$}
\end{figure}

In this section, we use the QP (\ref{qp.p1}, \ref{qp.p2}) away from the optimal solution as a means to determine the step in the shape normal direction and thus create an iterative solution technique very similar to SQP techniques known from linear spaces.
We solve the optimal control problem (\ref{qps.1}--\ref{qps.6}) by employing a CG--iteration for the reduced problem (\ref{oc.design}). I.e., we iterate over the variable $w$ and each time the CG--iteration needs a residual of equation (\ref{oc.design}) from $w^k$, we compute the state variable $z^k$ from (\ref{qps.2}--\ref{qps.6}) and then the adjoint variable $q^k$ from 
(\ref{ocad.1}, \ref{ocad.2}), which enables the evaluation of the residual
\begin{equation}
r^k:=-\left\llbracket f\right\rrbracket\left(p+\kappa p w^k+ q^k\right)+\mu\kappa-\mu\frac{\partial^2 w^k}{\partial\tau^2}\hspace{1mm}.
\end{equation}
The particular values for the parameters are chosen as $f_1=1000$, $f_2=1$ and $\mu=10$. The data $\bar{y}$ are generated from a solution of the state equation (\ref{oc2}, \ref{oc2}) with $u$ being the straight line connection of the points $(0.5,0)$ and $(0.5,1)$. 
The starting point of our iterations is described by a B--spline defined by the two control points $(0.6, 0.7)$ and $(0.4,0.3)$. We build a coarse unstructured tetrahedral grid $\Omega_h^1$ with roughly 6000 triangles as shown in the leftmost picture of figure \ref{fig-iterations}. We also perform computations on uniformly refined grids $\Omega_h^2,\Omega_h^3$ with roughly 24000 and 98000 triangles. In  figure \ref{fig-iterations} are also shown the next two iterations on the coarsest grid, where table \ref{fig-iterations} gives the distances of each shape to the solution approximated by
\[
\mathrm{dist}(u^k,u^*):= \int_{u^*}\left|\left\langle u^k,e_1\right\rangle-\frac12\right|ds
\]
where $u^*$ denotes the solution shape and $e_1=(1,0)$ is the first unit vector. Similar to \cite{VHS-shape-Riemann}, the retraction chosen for the shape is just the addition of the $q^kn_1$ to the current shape. In each iteration, the volume mesh is deformed according to the elasticity equation. Table \cite{VHS-shape-Riemann} demonstrates that indeed quadratic convergence can be observed on the finest mesh, but also that the mesh resolution has a strong influence on the convergence properties revealed.

The major advantage of the Newton method over a standard shape calculus steepest method based on the (reduced) shape derivative  
\[
dJ(y,u)[V]=-\int_u \left(\left\llbracket f\right\rrbracket p-\mu\kappa\right)\langle V,n\rangle ds
\]
is the natural scaling of the step, which is just 1 near to the solution. When first experimenting with a steepest descent method, we found by trial and error, that one needs a scaling  around 10 000 in order to obtain sufficient progress. 

\begin{table}[h]
\caption{Performance of shape Lagrange--Newton algorithms: distances $\mathrm{dist}(u^k,u^*)$ from the optimal solution on meshes with varying refinement. Quadratic convergence on the finest grid can be observed.}
\begin{center}
\begin{tabular}{r|c|c|c}
It.--No.&$\Omega_h^1$&$\Omega_h^2$&$\Omega_h^3$\\\hline\hline
0&0.0705945 & 0.070637 & 0.0706476\\
1&0.0043115 & 0.004104 & 0.0040465\\
2&0.0003941 & 0.000104 & 0.0000645\\
\end{tabular}
\end{center}
\end{table}

\section{Conclusions}
This paper presents a  generalization of the Riemannian shape calculus framework in \cite{VHS-shape-Riemann} to Lagrange--Newton approaches for PDE constrained shape optimization problems. It is based on the idea that Riemannian shape Hessians do not differ from classical shape Hessians in the solution of a shape optimization problem and that Newton methods still converge locally quadratically, if Hessian terms are neglected which are zero at the solution anyway. It is shown that this approach is viable and leads to computational methods with superior convergence properties, when compared to only linearly converging standard steepest descent methods. Nevertheless, several issues have to be addressed in future investigations, like:
\begin{itemize}
\item More refined retractions have to be developed for large shape deformations.
\item As observed during the computations, the shape deformation sometimes leads to shapes, where normal vectors can no longer be reliably evaluated. Provisions for those cases have be developed
\item Full Lagrange--Newton methods may turn out being not very computationally efficient. However, this paper lays the foundation for the construction of appropriate preconditoners for the reduced optimization problem in many practical cases.
\item The Riemannian shape space properties including quadratic convergence of the Lagrange--Newton approach seem to materialize only on very fine grids. A logical next development is then to use locally adapted meshes near the shape front to be optimized.
\end{itemize}
\section*{Acknowledgment} This research has been partly funded by the DFG within the collaborative project EXASOLVERS as part of the DFG priority program SPP 1648 SPPEXA. Furthermore, the authors are very grateful to an unknown referee for insightful comments, which helped significantly to improve the paper from an earlier version.
 
\bibliographystyle{plain}

\end{document}